%% file: deg17.092003.tex
\documentclass[a4paper]{article}

\newcommand{\septic}{septimic }
\newcommand{\septics}{septimics }

\usepackage[T1]{fontenc} 
\usepackage{theorem}
\newtheorem{prop}{\medskip{}{\bf Proposition}}[section]
\newtheorem{dfn}[prop]{\medskip{}{\bf Definition}}
\newtheorem{dfn-prop}[prop]{\medskip{}{\bf Definition-Proposition}}

\newtheorem{cor}[prop]{\medskip{}{\bf Corollary}}
\newtheorem{lemma}[prop]{\medskip{}{\bf Lemma}}

\usepackage{mathrsfs} 
\usepackage{rotating} 
\usepackage{psfrag} 
\usepackage{amssymb} 
\usepackage{multirow}

\newenvironment{preuve}{\medskip{}\noindent{\it Proof. }}{\hfill
    $\square$\medskip} \newenvironment{preuvesc}{\medskip{}\noindent{\it
      Proof. }}{\medskip}
\newenvironment{concsn}{\medskip{}\noindent{\bf Theorem }\it }{}
 \newcounter{fig}
\renewcommand{\caption}[1]{\addtocounter{fig}{1}
      \begin{center}
        {\sc Fig. \thefig} --- {\it #1}
\end{center}
}

\newcommand{\p}[2]{\mbox{$\mathbf
      P_{#1}^{#2}$}} \newcommand{\cp}[2]{\mbox{$\check{\mathbf
        P}_{#1}^{#2}$}} \newcommand{\sldfs}{\mbox{${\mathrm S\mathrm
        L}(2,{\mathbf F}_7)$}} \newcommand{\sldc}{\mbox{${\mathrm
        S\mathrm L}(2,{\mathbf C})$}} \newcommand{\sltc}{\mbox{${\mathrm
        S\mathrm L}(3,{\mathbf C})$}}
\newcommand{\K}[1]{\mbox{${\mathscr K}_{#1}$}}
\renewcommand{\H}[1]{\mbox{${\mathscr H}_{#1}$}}
\newcommand{\Kp}[1]{\mbox{$\K{}\,^\prime_{#1}$}}
\newcommand{\C}[1]{\mbox{${\mathscr C}_{#1}$}} \def\lr{\mathchar"2160}
\renewcommand{\l}[1]{\mbox{${\lr}_{#1}$}}

\newcommand{\I}{\mbox{$\mathscr I$}}
  \renewcommand{\O}{\mbox{$\mathscr O$}}
 \usepackage{xy}\xyoption{all}

\DeclareMathAlphabet{\mathbi}{T1}{cmr}{bx}{it}
\DeclareOldFontCommand{\bi}{\bfseries\itshape}{\mathbi}
\DeclareSymbolFont{letters}{OML}{ptmcm}{m}{it}
\DeclareMathSizes{12}{12}{7}{6}

\theorembodyfont{\rm} \newtheorem{rem}[prop]{\medskip{}{\bf Remark}}

\title{Degenerations of $(1,7)$--polarized abelian surfaces
}
\author{F. {\sc Melliez} and K. {\sc Ranestad}}
\date{}
\begin{document}
\maketitle

\section{Introduction}
The moduli space $A(1,7)$ of $(1,7)$--polarized abelian surfaces with
a level structure was shown by Manolache and Schreyer to be rational
with compactification $V(\K4)$ a Fano $3$-fold $V_{22}$ \cite{M-S}.
This compactification is obtained considering the embedding of the
abelian surfaces with their level structure in $\p{}{}V_{0}$, where
$V_{0}$ is the Schr\"odinger representation of the Heisenberg group of
level $7$.  Any such embedding is invariant under the action of
$G_{7}$, an extension of $H_{7}$ by an involution. The fixed points of
this involution and its conjugates in $G_{7}$ form an $H_{7}$-orbit of
planes $\p2+$ and $3$-spaces $\p3-$.  The $3$-fold $V(\K 4)$
parametrizes, what we denote by generalized $G_{7}$-embedded abelian
surfaces.  Every such surface intersects $\p2+$ in a finite subscheme
of length $6$, which we may classify by its type, namely the length of
its components.

In this paper we prove the

\medskip{}
\begin{concsn}\label{conc}
  Let $A$ a generalized $G_7$--embedded abelian surface of
  $\p{}{}V_{0}$, then according to the type of $\zeta_A=A\cap \p2+$
  the surface $A$ is

    \begin{center}
    \begin{tabular}{r|l}
      type of $\zeta_A$& description \\
      \hline
      $(1,1,1,1,1,1)$ & smooth and abelian \\
      $(2,1,1,1,1)$ & translation scroll $(E,\pm\sigma)$ with $2\cdot
\sigma \neq 0$\\
      $(3,1,1,1)$ & tangent scroll $(E,0)$\\
      $(2,2,2)$ & doubled translation scroll $(E,\pm\sigma)$ with
$2\cdot \sigma=0$ and $\sigma \neq 0$\\
      $(2,2,1,1)$ & union of seven quadrics\\
$(4,2)$ & union of seven doubled projective planes\\
$(2,2,2)_s$ & union of fourteen projective planes
    \end{tabular}
\end{center}
\end{concsn}
The two occurring $(2,2,2)$ cases (abusively denoted by $(2,2,2)$ and
$(2,2,2)_s$) are special ones and are described in {\sc Fig.} 1
(subsection \ref{degen}).

\medskip{}

In the first section we first recall some basic facts on
$(1,7)$--polarized abelian surfaces with level structure, and
construct a compactification of their moduli space which turns out to
be the same as the one constructed by Manolache and Schreyer. This
alternative proof of their result (note that the main argument appears
in their paper) is inspired by the case of $(1,5)$ polarizations
\cite{BHM}, in particular we do not use syzygies (although both
constructions coincide in the $(1,5)$ case!).

The second section deals essentially with remarks on prime Fano
threefolds of index $12$. For the one which is endowed by a faithful
action of $\p{}{}\sldfs$ we give a careful description of the
boundary.

In the last section we prove the theorem by interpreting the boundary
in terms of surfaces in $\p{}6=\p{}{}V_0$.

Note that A. Marini investigates also such degenerations in \cite{AMa}
from a different interpretation of $V(\K4)$, namely that of twisted
cubics apolar to the ``Kleinian'' net of quadrics.


\subsection*{Notations}

\indent The base field is the one of complex numbers ${\bf C}$.  If
$R$ is a vector space, the Veronese map from $R$ to $S^nR$ (as well as
its projectivisation) will be denoted by $\nu_n$ :
$$\xymatrix{R\ar[r]^{\nu_n}&S^nR}.$$
If $s\in {\rm Hilb} (n,{\bf P}R)$
the type of $s$ ({\it i.e.} the associated length partition of $n$)
will be labelled $\lambda_s$ :
$$\xymatrix{s\ar[r]^{\lambda}&\lambda_s}.$$
If $H$ is a hypersurface
of ${\bf P}R$ then $e_H=0$ is an equation of $H$.

The irreducible representations of $\sldfs$ will be denoted by ${\bf
  C}$, $W_3$, $W_3^\lor$, $U_4$, $U_4^\lor$, $W_6$, $U_6$, $U_6^\lor$,
$W_7$, $W_8$ and $U_8$. The algebra of representations of the group
$\sldfs$ is a quotient of $${\bf Z}[{\bf
  C},W_3,W_3^\lor,U_4,U_4^\lor,W_6,U_6, U_6^\lor,W_7,W_8,U_8]$$
where
$\bf C$ denotes the trivial representation, $W_n$ denotes an
irreducible $\p{}{}\sldfs$--module of dimension $n$ and $U_n$ denotes
an irreducible $\sldfs$--module of dimension $n$ on which $\sldfs$
acts {\it faithfully}.

The corresponding table of multiplication can be found in \cite{M-S}
and \cite{Dol} with the following possible identifications

\begin{center}
\begin{tabular}{|c|c|c|c|c|c|c|c|c|c|c|c|}
\hline
\cite{Dol} &$V_1$&
$V_3=V_{-}$&$V_3^*$&$V_4=V_{+}$&$V_4^*$&$V_6$&$V'_6$&$V'^{*}_6$&$V_7$&
$V'_8$&$V_8$\\
\hline
\cite{M-S} &I  &$W$&$W'$&$U$&$U'$&$T$&$T_1$&$T_2$&$L$&$M_1$&$M_2$\\
\hline
$\times$   &${\bf C}$&
$W_3$&${W_3^{\lor}}$&$U_4$&$U_4^{\lor}$&$W_6$&$U_6$&$U_6^{\lor}$&$W_7$
&
$U_8$&$W_8$\\
\hline
\phantom{\Large P}& $\cdot$ &
\p2+&\cp2+&\p3-&\cp3-&\p5+&\p5-&\cp5-&\p6+& \p7-&\p7+\\
\hline
\end{tabular}
\end{center}

Note that what are denoted by $\p2+$ and $\p3-$ are respectively
denoted by $\p-2$ and $\p+3$ in \cite{GPo}.

The seven dimensional vector spaces $V_0,...,V_5$ are irreducible
$H_7$--modules, a useful multiplication table have been stolen (among
many other things) in \cite{M-S}.
\section{Moduli space : a compactification}

\bigskip Let $A$ be an abelian surface, {\it i.e.} a projective
complex torus ${\bf C}^2/\Lambda$ where $\Lambda$ is a (maximal)
lattice of ${\bf C}^2\simeq {\bf R}^4$.  Then the variety ${Pic}^0(A)$
turns to be an abelian surface as well (isomorphic to $({\bf
  C}^2)^\lor / \Lambda^\lor$); this latter one is called the
\textit{dual} abelian surface of $A$ and will be denoted by $A^\lor$.
As additive group, the surface $A$ acts on itself by translation, if
$x\in A$ we will denote by $\tau_x$ the corresponding translation.

A line bundle of type $(1,7)$ on $A$ is the data of an ample line
bundle $\mathscr L$ such that the kernel of the isogeny
$$\xymatrix{\varphi_{\mathscr L} : A\ar[r]&A^\lor, &
  x\ar[r]&\tau_x^*\mathscr L\otimes\mathscr L^{-1} }$$
is the affine
plane (with origin) over $\mathbf F_7$, {\it i.e.} such that
$ker(\varphi_{\mathscr L})\simeq \mathbf Z_7 \times \mathbf Z_7$.

A $(1,7)$--polarization on $A$ is an element of

$$\{(A,\varphi_{\mathscr L})\ |\ \mathscr L\ \mbox{\rm is of type}\ 
(1,7) \}.$$

Thanks to Mumford, a coarse moduli space of $(1,7)$--polarized abelian
surfaces exists, we will denote it by $M(1,7)$.

Now choose a {\it generic} abelian surface, say $A$, then
$V_0={\mathrm H}^0(A,\mathscr L)$ is of dimension $7$.
The group $ker(\varphi_{\mathscr L})\simeq \mathbf Z_7 \times \mathbf
Z_7$ becomes a subgroup of $\p{}{}\mathrm{SL}(V_0)$.
It is certainly safer to work with linear representations rather than
projective ones so we need to lift the action of $\mathbf Z_7 \times
\mathbf Z_7$ on $\p{}{}V_0$ to an action of one of its central
extensions on $V_0$.  The Schur multiplier of $\mathbf Z_7 \times
\mathbf Z_7$ is known to be $\mu_7$ so any projective representation
of $\mathbf Z_7 \times \mathbf Z_7$ is induced by a linear
representation of what is called the ``Heisenberg group of level $7$''
and denoted by $H_7$: that is to say for all $n\in \bf N^*$ and all
projective representation $\rho$ we get a cartesian diagram:

$$\xymatrix@1@W=1cm@H=.6cm{ \ar[r]1&{\mu_n}\ar[r]& {\rm SL}(n,\bf
  C)\ar[r] & {\rm
    \p{}{}SL}(n,\bf C)\ar[r]& 1 \\
  1\ar[r]&{\mu_7}\ar[r]\ar[u]&H_7\ar[r]\ar[u]&\ar[r] \mathbf Z_7
  \times \mathbf Z_7\ar[u]_{\textstyle{\rho}}&1}$$

In this way $V_0$ becomes a $H_7$-module (of rank $7$), this
representation is called the ``Schr{\"o}dinger'' representation of
$H_7$. We have now a way to identified {\it all} the vector spaces
${\mathrm H}^0(A',\mathscr L)$ for any abelian surface $A'\in M(1,7)$
as they are all isomorphic to $V_0$ as $H_7$-modules.  This looks too
good to be true.
So what is wrong? We implicitly made an identification between
$ker(\varphi_{\mathscr L})$ and $\mathbf Z_7 \times \mathbf Z_7$ and
this is certainly defined up to \sldfs\ only!  So the construction is
only invariant under $N_7=H_7\rtimes \sldfs$ which turns out to be the
normaliser of $H_7$ in $\mathrm{SL}_7(\bf C)\simeq\mathrm{SL}(V_0)$
(this seems to be a collective agreement...).

So to any basis $s$ of $ker(\varphi_{\mathscr L})$ corresponds an
embedding
$$\zeta_s:A\longrightarrow \p{}{}V_0.$$
The group $\sldfs$ acts on the
set of basis of $ker(\varphi_{\mathscr L})$ and we immediately get
another complication (which will turn out to be quite nice after all):
$$\zeta_s(A)=\zeta_{-s}(A).$$
Let us denote by $G_7=H_7\rtimes
\{-1,1\}\subset N_7$. This group (after killing $\mu_7$) is in general
the full group of automorphisms of the surface $\zeta_s(A)$: if $b$ is
any element of $\mathbf Z_7 \times \mathbf Z_7$ and
$\tau_b:\p{}{}V_0\longrightarrow \p{}{}V_0$ is the involution induced
by the corresponding ``${-1}$'' of $G_7$, then $\tau_b$ leaves
$\zeta_s(A)$ (globally) invariant and is induced by the ``opposite''
map $x\mapsto -x$ on $A$ for a good choice of the image of the origin
on $\zeta_s(A)$. In other words, $\tau_b \cdot \zeta_s = \zeta_{-s}$.

As the cardinality of $\sldfs/\{-1,1\}=\p{}{}\sldfs$ is 168, each
element of $M(1,7)$ will be embedded in $\p{}{}V_0$ in $168$ ways
(distinct in general). We get a brand new moduli space considering a
$(1,7)$--polarized abelian surface together with one of its embedding,
this moduli space will be denoted by $A(1,7)$:
$$A(1,7)=\{((A,\varphi_{\mathscr L}),s) |\, (A,\varphi_{\mathscr
  L})\in M(1,7),\ s \ \mbox{is a basis of} \ ker(\varphi_{\mathscr L})
\}/\square$$
in which the equivalence relation $\square$ is the
expected one, $(X_1,s_1)\ \square\ (X_2,s_2)$ if
$\zeta_{s_1}(X_1)=\zeta_{s_2}(X_2)$ (fortunately, this implies
$X_1=X_2$).  The choice of basis (or embedding) is the level structure
refered to in the introduction.

\renewcommand\theenumi{{\it \roman{enumi}}} \bigskip

Here are useful remarks:
\begin{enumerate}\label{useful}
\item The surface $\zeta_{s}(A)$ is of degree 14;
\item by construction if $x\in A$ then the set of 49 points
  $\zeta_s(\varphi_{\mathscr L}^{-1}(\varphi_{\mathscr L}(x)))$ is an
  orbit under the action of $H_7$ (or $H_7/\mu_7$ if we want to be
  precise);
\item the above construction works as well for elliptic curves, so in
  particular $\p{}{}V_0$ contains naturally $G_7$--invariant embedded
  elliptic curves (of degree 7);
\item if $b\in\mathbf Z_7 \times \mathbf Z_7$ the involution $\tau_b$
  induces a $\sldfs$--module structure on $V_0$, as such a module
  $V_0$ splits in $V_0=W_3\oplus U_4$ where both $W_3$ and $U_4$ are
  irreducible $\sldfs$--modules of dimension $3$ and $4$ respectively
  (such that $S^3W_3^\lor\simeq S^2U_4$). The projective plane
  $\p{}{}W_3$ and the projective space $\p{}{}U_4$ in $\p{}{}V_0$ are
  point by point invariant by the involution $\tau_b$.  For a given
  $b\in\mathbf Z_7 \times \mathbf Z_7$ these two spaces will often be
  denoted by $\p2+$ and $\p3-$ (the signs come from the following:
  $W_3$ is also a $\p{}{}\sldfs$--module, {\it i.e.} $-1\in \sldfs$
  acts trivially on it, but $\sldfs$ acts faithfully on $U_4$);
\item if $E$ is a $G_7$--invariant elliptic curve in $\p{}{}V_0$ then
  the curve $E$ intersects any $\p2+$ in one point (corresponding to
  the image of $0$) and any $\p3-$ in $3$ points (corresponding to its
  non trivial $2$--torsion points);
\item the latter holds also for abelian surfaces, with decomposition
  $6+10$ corresponding to the odd and even $2$--torsion points
  (\cite{LB});
\item by adding a finite set of $G_{7}$--invariant heptagons to the
  union of the $G_7$--invariant embedded elliptic curves of degree
  $7$, one gets a birational model of the Shioda modular surface of
  level $7$.  It intersects each $\p2+$ in a plane quartic curve
  $\Kp4$, the socalled Klein quartic curve (\cite{Mie} or \cite{Fr}
  which contains original references to Klein).
\end{enumerate}

Following what happens in the $(1,5)$ case we next consider the
rational map
$$\xymatrix{\kappa:\p{}{}V_0\ar@{-->}[r] &
  \p{}{}(\mathrm{H}^0({\mathscr O}_{{\bf P}V_O}(7))^{G_7})^\lor}$$
{\it i.e.} the blowup of $\p{}{}V_0$ by the linear system of
$G_7$--invariant \septics. {\it In what follows, by a
  `$G_7$--invariant \septic' we always (and imprudently) mean a
  \septic in this linear system}. The isomorphism of
$\p{}{}\sldfs$--modules $\mathrm{H}^0({\mathscr O}_{{\bf
    P}V_O}(7))^{G_7}\simeq W_7\oplus {\bf C}$ shows that $\kappa$
takes, {\it a priori}, its values in a $\p{}7$. We will see that the
image of $\kappa$ is in fact contained in $\p{}{}W_7$.  First we
analyse the base locus of these \septic hypersurfaces.

\begin{lemma}
  A $G_7$--invariant \septic hypersurfaces of $\p{}{}V_0$ contains any
  of the forty--nine projective spaces \p3-.
\end{lemma}
\begin{preuve}
  Obviously the vector space $\mathrm{H}^0({\mathscr O}_{{\bf
      P}V_O}(7))^{G_7}$ is a $\p{}{}\sldfs$--module, so considering
  the restriction to any projective space $\p3-=\p{}{}U_4$ we get a
  map $\mathrm{H}^0({\mathscr O}_{{\bf P}V_O}(7))^{G_7}\longrightarrow
  \mathrm{H}^0({\mathscr O}_{{\bf P}U_4}(7))=S^7U_4^\lor$ which need
  to be $\sldfs$--equivariant (the entire collection of $\p3-$'s being
  invariant under the action of $G_7$).  But $U_4$ is a faithful
  module for $\sldfs$ and $7$ is odd,so the map is the zero map.
\end{preuve}

Using Bezout theorem we get
\begin{cor} A $G_7$--invariant \septic hypersurface of
  $\p{}{}V_0$ contains any $G_7$--invariant elliptic curve of
  $\p{}{}V_0$ as well as its translation scroll by a non trivial
  $2$--torsion point.
\end{cor}

Notice that our forty nine $\p2+$ constitute an orbit under $G_7$, so
it makes sense to consider \textit{the} surface $\kappa(\p2+)$. As any
objects we are interested in intersects these planes, this will
certainly be a good point to understand this surface.
\begin{cor}
  `The' plane $\p2+$ is mapped by $\kappa$ to a Veronese surface of
  degree nine in $\p{}{}W_7$.
\end{cor}
\begin{preuve}
  Let $P\in{\mathrm H}^0({\mathscr O}_{{\bf P}V_0}(7))^{G_7}$. From
  the preceding corollary the hypersurface $\{P=0\}$ contains the
  Shioda modular surface, and in particular its $7$--torsion points
  sections. So the \septic hypersurface $\{P=0\}$ intersects any plane
  $\p2+$ along a reducible \septic curve, union of the Klein quartic
  curve (which is the curve of $7$--torsion points contained in a
  plane $\p2+$) and a cubic curve.  The irreducible
  $\p{}{}\sldfs$--modules decomposition of the vector space
  $S^3W_3^\lor\simeq W_7\oplus W_3$ brings us to the proclaimed
  result, if one notices the equality $ W_7\simeq W_7^\lor$.
\end{preuve}

\begin{rem}
  This phenomenon holds also in the $(1,5)$ case where $\p2+$ is
  mapped to a (projected) Veronese surface of degree $25$ in a
  Grassmannian ${\rm Gr}(1,\p3{})\subset\p5{}$ known as the bisecants
  variety of a certain rational sextic curve in $\p3{}$. The image of
  the blow-up of `the' line $\p1-$ intersects the latter Grassmannian
  along the sextic complex of lines contained in a plane of a {\it
    dual} sextic of $\p3{}$. In this case, any $G_5$--embedded
  $(1,5)$--polarized abelian surface is mapped to a ten--secant plane
  to $\kappa (\p2+)$ (which intersects the sextic complex along six
  lines).  Although the same kind of results are expected in our
  situation, here is a difference between the two cases, namely
\end{rem}

\begin{rem}\label{smoore}
  The vector space $\mathrm{H}^0({\mathscr O}_{{\bf P}V_O}(7))^{G_7}$
  is not given by determinants of (symmetric) Moore matrices
  (\cite{GPo}).
\end{rem}
\begin{preuve}
  Let us recall first what is a Moore matrix; we have a nice
  isomorphism of irreducible $N_7$--modules (defined up to homothety)
  $S^2 V_4= U_4\otimes V_0$ which induces a map $U_4 \longrightarrow
  S^2V_4\otimes V_0^\lor.$ For a good choice of basis in $V_0^\lor$ we
  get a $7\times 7$ matrix with coefficients in $V_0^\lor$ which is
  called a (symmetric) ``Moore matrix''.  Now considering determinants
  ({\it i.e.} the map $\xymatrix{{S^2V_4}\ar[r]^{\ `` S^2 \Lambda^7
      "}& {\bf C}}$) we get a first map
  $$S^7U_4\longrightarrow S^7V_0^\lor={\mathrm H}^0({\mathscr O}_{{\bf
      P}V_0}(7)),$$
  which composed with the projection to the
  invariant part yields a map $$S^7U_4\longrightarrow {\mathrm
    H}^0({\mathscr O}_{{\bf P}V_0}(7))^{G_7}.$$
  This latter one is
  certainly zero: the action of $-1\in\sldfs$ cannot be trivial on any
  $\sldfs$--invariant subspace of the vector space $S^7U_4$.
\end{preuve}


Nevertheless {\it anti-symmetric} Moore matrices play a fundamental
role in the $(1,7)$ case.  They are defined by the isomorphism of
irreducible $N_7$--modules $\Lambda^2 V_4= W_3^\lor\otimes V_0$. The
locus (in $\p{}{}V_0$) where such matrix drops its rank are
Calabi-Yau threefolds (see \cite{GPo}) and will appear in subsection
\ref{ANTI}.

\bigskip

\begin{prop}\label{suite} 
  The image by the map $\kappa$ of a $G_7$--embedded
  $(1,7)$--polarized abelian surface is (generically) a projective
  plane and we have a factorization
$$\xymatrix{{A^\flat}\ar[r]^{49:1\ }&A^{\lor\flat} \ar[r]^{2:1\ }
  &K_{A^\lor}^\flat\ar[r]^{2:1\ }&\kappa (A)}$$
where the surface
${A^\flat}$ is the blowup of the surface $A$ along its intersection
with the base locus of the $G_7$-invariant \septics.
\end{prop}
\begin{preuve}
  Assume the $(1,7)$--polarization of $A\in M(1,7)$ is given by a {\it
    very} ample line bundle, then from the $G_7$--equivariant
  resolution of the surface $A$ in $\p{}{}V_0$ which can be found in
  \cite{M-S}, one can check (or read, it is easier {\it op. cit.},
  appendix) that ${\rm dim}({\rm H^0}({\mathscr O}_A(7))^{G_7})=3$.

Supposing the map $\kappa|_A$ is finite (as it should be) then we have a factorization
$$\xymatrix{{A^\flat}\ar[r]^{49:1\ }&A^{\lor\flat} \ar[r]^{2:1\ }
  &K_{A^\lor}^\flat\ar[r]^{2:1\ }&\kappa (A)}.$$
The first two maps
(as well as their degree) come from the construction itself, the
degree of the last one is due to {\'E}tienne Bezout.

If $\kappa|_A$ is not finite then it is composed with a pencil.  We
may assume that ${\rm Pic} (A)$ has rank $1$, \textit{i.e.} all curves
are hypersurface sections or translates thereof.  But no such curve is
fixed by $G_7$, so $\kappa|_A$ has at most isolated base points.  So
the linear system is a subsystem of $|7\cdot h|$.  In particular if it
composed with a pencil the fibers are curves in $|h|$ through the base
locus.  But there are clearly no such curve, since already the $49$
$\p3-$ span $\p{}6=\p{}{}V_0$.  Thus the map $\kappa$ is finite on $A$
and the proposition follows.
\end{preuve}

Let us denote by ${A(1,7)^v}$ the (open) subset of ${A(1,7)}$
corresponding to $(1,7)$--polarized abelian surfaces for which the
polarization is given by a very ample line bundle. Using the preceding
section we get a nice embedding ${A(1,7)^v}$ in the (compact) variety
of six secant planes to $\kappa (\p2+)$ and then a natural
compactification $\overline{A(1,7)^v}$ (see N. Manolache and F.-O.
Schreyer's paper \cite{M-S} for another proof) which comes with a
geometric meaning; we just need to consider the proper transform of a
six--secant plane of the Veronese surface by $\kappa^{-1}$. Moreover,
{\it any} $(1,7)$--polarized abelian surface is mapped into the
hyperplane $\p{}{}W_7$ of $\p{}{}{\mathrm H}^0({\mathscr O}_{{\bf
    P}V_0}(7))^{G_7}$ so their union is contained in a \septic
hypersurface of $\p{}{}V_0$. So we have

\begin{cor}\label{1.8}
  The moduli space $\overline{A(1,7)^v}$ is isomorphic to the unique
  prime Fano threefold of genus $12$ which admits $\p{}{}\sldfs$ as
  its automorphisms group. The universal $(1,7)$--polarized abelian
  surface with level $7$ structure is birational to the unique
  $N_7$--invariant \septic hypersurface of $\p{}{}V_0$.
\end{cor}
\begin{preuve} Let us denote by $X_7$ the unique
  $N_7$--invariant \septic{} hypersurface of $\p{}{}V_0$ and by
  $B_\kappa$ the base locus of the $G_7$--invariant septimic
  hypersurfaces. Put $Y_4=\overline{\kappa (X_7\backslash
    B_\kappa)}\subset \p{}{}W_7$ and consider the diagram
  \label{diagr}
  $$\xymatrix{&&I\ar[rd]^{p_2}\ar[ld]_{p_1}\\
    X_7\backslash B_\kappa\ar[r]^\kappa &Y_4&&{\bf G} (3,\, W_7)}$$
  where $I\subset \p{}{}W_7\times {\bf G} (3,\, W_7)$ denotes the
  graph of the incidence correspondence between $\p{}{}W_7$ and the
  (projective) fibers of the tautological sheaf over ${\bf G} (3,\,
  W_7)$ and where $p_1$ and $p_2$ are the natural projections.  In
  order to prove birationality we just need to prove that a general
  point of $X_7$ is contained in one (and only one) abelian surface.
  One first needs to remark, using representation theory for instance,
  that both the hypersurfaces $X_7$ and $Y_4$ are irreducible.
  
  Let $A\in {A(1,7)^v}$ a $G_7$--embedded abelian surface. We have:
 \begin{itemize}
 \item the septimic $X_7$ contains the surface $A$;
 \item the surface $A$ intersects $\p2+$ along a reduced scheme;
 \item the surface $A$ {\it is not contained} in the base locus
   $B_\kappa$.
\end{itemize}

The only non obvious fact is the third item. But $B_\kappa$ intersects
$\p2+$ along a Klein quartic curve $\Kp4$ so if we had $A\subset
B_\kappa$ this would imply the non emptiness of $A\cap \Kp4$ and in
such cases $A\cap \p2+$ admits a double point (see e.g. section
\ref{ftv} below) contradicting the second item.  Next the map $A
\longmapsto A\cap\p2+$ is injective (see \cite{M-S}) so the plane
$\overline{\kappa (A \backslash B_\kappa)}$ entirely characterizes the
surface $A$.  Summing up all we know, we can pretend that two
distincts surfaces $A$ and $A'$ intersect each other either on
\begin{itemize}
\item the threefold $B_\kappa$ (which is of codimension $2$ in $X_7$);
\item the preimage by $\kappa$ of the points in $Y_4\subset \p{}{}W_7$
  which are contained in more than one six-secant plane to the
  Veronese surface $\overline{\kappa (\p2+\backslash \Kp4)}$.
\end{itemize}

Since $A$ is not contained in $B_\kappa$, it remains to show that $A$
is not contained in the second locus. But one proves easily that
the second locus is $2$-dimensional, being the preimage of the reunion
of the Veronese surface $\overline{\kappa (\p2+\backslash \Kp4)}$
itself and its ruled surface of trisecant lines (for which the base is
isomorphic to the Klein quartic curve $\K4$ of the dual plane
$\cp2+$).
\end{preuve}

\begin{rem}
  Notice that one can also prove (using Schubert calculus) that the
hypersurface $Y_4$ has degree four in $\p{}{}W_7$ (this is true for
any collection of six-secant planes to such projected Verenose
surface).
\end{rem}

\section{Fano threefolds $V_{22}$}\label{ftv}
Let us recall a characterization of prime Fano threefolds of genus 12
(cf \cite{Mu}).

\begin{dfn-prop}
  Any Fano threefold of index 1 and genus 12 is isomorphic to the
  variety of sums of power $VSP
  (F,6)=\overline{\{({\l1},...,{\l{6}})\in {\rm Hilb}_6\, \p{}{}W^\lor
    \ | \ 
    e_F\equiv\mathop{\Sigma}\nolimits_{i=1}^{6}e_{{\lr}_{i}}^4\}}$ of
  a plane quartic curve $F$.  Conversely, if $F$ is not a Clebsch
  quartic (i.e. its catalecticant invariant vanishes), then $VSP
  (F,6)$ is a Fano threefold of index 1 and genus 12.  Its
  anti-canonical model is denoted by $V_{22}$.
\end{dfn-prop}

\subsection{Construction}
Let $W$ be an irreducible $\sltc$-module of dimension $3$, we have a
decomposition of $\sltc$-modules (\cite{F-H})
$$S^2(S^2W)^\lor=S^4W^\lor\oplus S^2W$$
generating an exact sequence
$$0\longrightarrow S^4W^\lor\longrightarrow {\rm
  Hom}(S^2W,S^2W^\lor).$$
So a plane quartic $F$ in $\p{}{}W$ whose
equation is given by `an' element $e_F$ of $S^4W^\lor$ gives rise to
`a' morphism $\alpha_F:S^2W \longrightarrow S^2W^\lor$ and a quadric
$\mathbf Q_F$ in $\p{}{}(S^2W)$ of equation $\alpha_F(x)\cdot x=0$ or
even $x\cdot \alpha_F(x)=0$ by the canonical identification
$S^2W=(S^2W^\lor)^\lor$. From the equality ${\rm h}^0({\mathscr
  I}_{\nu_2(F)}(2))=7$, we get a \textit{characterization of this
  quadric} by the two properties: \renewcommand\theenumi{{\it
    \roman{enumi}}}
\begin{enumerate}
\item the two forms on $W$ defined by $\alpha_F(\nu_2(-))\cdot
  \nu_2(-)$ and $e_F$ are proportional {\it i.e.} the quadric $\mathbf
  Q_F$ and the Veronese surface $\nu_2(\p{}{}W)$ intersects along the
  image of the plane quartic $F$ under $\nu_2$;
\item the quadric $\mathbf Q_F$ is apolar to the Veronese surface
  $\nu_2(\p{}{}W^\lor)$ of $\p{}{}S^2W^\lor$ {\it i.e.} apolar to each
  element of the vector space ${\rm H}^0({\mathscr I}_{\nu_2({\bf
      P}W^\lor)}(2))\simeq S^2W^\lor\subset S^2(S^2W)$.
\end{enumerate}

\begin{lemma}[Sylvester]
  The minimal integer $n$ for which $VSP(F,n)$ is non empty is the
  rank of $\alpha_F$ (called the catalecticant invariant of the
  quartic curve).
\end{lemma}

\begin{preuve}
  This well known result of Sylvester (see e.g. Dolgachev and Kanev
  \cite{D-K}, Elliot \cite[page 294]{Ell}) can be deduced from the
  following observation: let $n\in {\bf N}^*$, then
  
  $$\nu_2(VSP(F,n))=\{({p_1},...,{p_n})\in VSP({\mathbf Q}_F,n) \ | \ 
  p_{\times} \in \nu_2(\p{}{}W^\lor)\}.$$
  Indeed if
  $\square=({\l1},...,{\l{n}})\in VSP(F,n)$ then
  $e_F\equiv\mathop{\Sigma}\limits_{i=1}^{n}e_{\l{i}}^4$ for a good
  normalization of $e_{l_\times}$ and the quadric ${\bf Q}\subset
  \p{}{}S^2W$ of equation $$e_{\bf
    Q}\equiv\mathop{\Sigma}\limits_{i=1}^{n}e_{\nu_2(\lr_{i})}^2$$
  is
  endowed with the two properties which characterize the quadric ${\bf
    Q_F}$: the second one is a direct consequence of ${\rm
    H}^0({\mathscr I}_{\nu_2({\bf P}W^\lor)}(2))\subset{\rm
    H}^0({\mathscr I_{\nu_2(\square)}}(2))$ and the first one arises
  by construction. Applying $\nu_2^{-1}$ we get the required equality.
\end{preuve}

Define the vector space $Y_\lr\subset S^2W$ such that the line $\lr$
of the plane $\p{}{}W$ induces the exact sequence
$$\xymatrix{0\ar[r]&\lower2mm\hbox{${\bf C}\cdot e_{\lr}^2$}\ar[r]
  &S^2W^\lor\ar[r]&\lower3mm\hbox{$Y_{\lr}^\lor$}\ar[r]&0},$$
that is
to say $Y_\lr$ is the orthogonal space (in $S^2W$) of $e_{\lr}^2$.

\begin{dfn}
  The subscheme ${\mathscr C}_{\lr}$ of the plane $\p{}{}W^\lor$
  defined by ${\mathscr C}_{\lr}=\{x\in \p{}{}W^\lor \ | \ 
  e_x^2\in\alpha_F(Y_{\lr})\}=\nu_2^{-1} (\alpha_F (Y_\lr))$ is called
  the {\it anti--polar conic} of the line $\lr$ (with respect to the
  quartic $F$).
\end{dfn}

\noindent Alternatively, if $\alpha_F$ has maximal rank we have obviously ${\mathscr
  C}_{\lr}=\{x\in \p{}{}W^\lor \ | \ e_x^2\cdot \alpha_F^{-1}
(e_\lr^2) =0\}$.

\noindent Set $n={\rm rank}(\alpha_F)$; the construction of a point of
$VSP(F,n)$ is now very easy by the following

\begin{cor}\label{conics}
  A point $(\l{1},...,\l{n})$ lies inside $VSP(F,n)$ if and only if
  $\l{i}\in\C{\lr_{j}}$ when $i\neq j$,
\end{cor}
which is an easy consequence of the classical construction of a point
of $VSP({\bf Q},n)$ when $\bf Q$ is a quadric of rank $n$. 

\medskip

We will need to understand

\subsection{Conics on the anti-canonical model}

Let $V$ be the seven dimensional vector space defined by the exact
sequence $\xymatrix{0\ar[r]&W\ar[r] &S^3W^\lor\ar[r]^{\ \ 
    p_F}&V\ar[r]&0}$ where the second map is induced by $F\in
S^4W^\lor\subset {\rm Hom} (W,S^3W^\lor)$ and denote by ${\mathscr
  V}_{2,9}'$ the image of $\p{}{}W^\lor$ in $\p{}{}V$ by the Veronese
embedding $\nu_3$ composed with the third map $p_F$.

By definition $\square=(\l{1},...,\l{n})\in VSP(F,6)$ if and only if
the image by $p_F$ of the six dimensional vector space (in
$S^3W^\lor$) spanned by $e_{\lr_{i}}^3$ is of rank $3$. Thus we get a map of $VSP(F,6 )$ into the Grassmannian $G (3,V)$, by
$(\l1,..,\l6)\mapsto p_F(<\l1,...,\l6>)$.

\begin{rem}
  { The image of $VSP(F,6)$ in the Pl\"ucker embedding of the Grassmannian is the anti-canonical model  $V_{22}$ of this Fano
    threefold, it is isomorphic to the variety of $6$ secant planes to
    the projected Veronese surface ${\mathscr V}_{2,9}'$.}
\end{rem}

Now it is reasonable to talk about conics on $V_{22}$.

Denote by $F^\flat$ the dual quartic of $\p{}{}W$ of equation
$\alpha_F^{-1}(e_\lr^2)\cdot e_\lr^2=0$, in other words we have
$F^\flat=\{\lr\in\p{}{}W^\lor\, |\, \lr\in \C{\lr} \}$ (the quartic
$F^\flat$ reduces to a doubled conic when $n=5$) and by $H_F$ the
sextic of $\p{}{}W^\lor$ given by $H_F=\{\lr\in\p{}{}W^\lor\, |\, {\rm
  rank} (\alpha_F ^{-1}(e_\lr^2))\leqslant 1 \}$.  Let
$\lr\in\p{}{}W^\lor\backslash\ H_F$, then the anti--polar conic
$\C{\lr}$ is smooth and we can write the abstract rational curve
$\overline{\C{\lr}}$ as $\p{}1=\p{}{}S_1$ with $\dim S_1=2$, {\it
  i.e.} we put $S_1={\rm H}^0 (\O_{\overline{\mathscr C}_{\lr}} (1))$. Put
$S_n:=S^nS_1$, then $\p{}{}S_n$ is identified with the divisors of
degree $n$ on $\p{}{}S_1$ and we have
\begin{lemma}\label{div5}
  The set of divisors of degree $5$ on the anti--polar conic $\C{\lr}$
  given by \hbox{$\{\lr'+\C{\lr}\cap \C{\lr'},\, {\lr'\in\C{\lr}}\}$}
  is a projective line in $\p{}{}S_5$. This $g_1^5$ admits a base
  point if $\lr\in F^\flat$.
\end{lemma}

\begin{preuve}
  Let $D$ be such a divisor. By corollary \ref{conics} the divisor $D$
  is completely determined by {\it any one} of its (sub)--divisor of
  degree $1$, so the variety of such divisors is a curve of first
  degree in $\p{}{}S_5$.
\end{preuve}

\begin{cor}
  The points of $VSP (F,6)$ which contain a given line $\lr$ describe
  a conic ${\mathcal C_\lr}$ on the anti-canonical model $V_{22}$. The
  two conics ${\mathcal C_\lr}$ and $\C{\lr}$ have the same rank.
\end{cor}

\begin{preuve}
  Let $\lr$ outside $H_F$. Then the image of $\C{\lr}$ by $\nu_3$ is a
  rational normal sextic projected by the map $p_F$ to a smooth sextic
  inside $\p{\lr}{4}:= \p{}{} (p_F ({\rm H}^0 (\O_{\mathscr C_{\lr}}
  (6))))$ --- we have an injection ${\rm H}^0 (\O_{\mathscr C_{\lr}}
  (6)) \subset S^3W^\lor$ and it is a simple matter to check ${\rm
    H}^0 (\O_{\mathscr C_{\lr}} (6))\cap \ker (p_F)$ is of dimension
  $2$, moreover identifying $\ker (p_F)$ and $W$ we have $\p{}{} ({\rm
    H}^0 (\O_{(\mathscr C_{\lr})} (6))\cap W)=\lr\subset \p{}{}W$ ---.
  Now $\p{\lr}{4}$ also contains the image of $\nu_3 (\lr)$ and
  projecting from this latter point the sextic becomes :
    \begin{enumerate}
    \item a rational sextic on a quadric of a $\p{}3$ generically;
    \item a rational quintic on a quadric of a $\p{}3$ if $\lr\in
      F^\flat$.
    \end{enumerate}
    These curves are obviously on a quadric, since a six secant plane to
    ${\mathscr V}_{2,9}'$ passing through the point $p_F (\nu_3
    (\lr))$ will be mapped to a 5-secant ({\it resp} 4-secant) line to
    this rational curve.
    
    Now if $\lr\in H_F$, $\C{\lr}$ breaks in two lines, say $\l1$ and
    $\l2$. We get two systems of six secant planes to ${\mathscr
      V}_{2,9}'$ containing $p_F (\nu_3 (\lr))$, one of these
    intersects $p_F (\nu_3 (\l1))$ in two {\it fixed} points and
    intersects the twisted cubic $p_F (\nu_3 (\l2))$ along a pencil of
    divisors of degree $3$. In particular, such collection is mapped
    to a line by $\kappa$.
\end{preuve}

\begin{cor}\label{CY}
  If $a$ is not on $H_F$ the threefold $p_1p_2^{-1} ({\mathcal C_a})$
  is a quadric cone $\Gamma_a$ of $\p{a}{4}$. If $a\in H_F$ the cone
  $\Gamma_a$ splits in two $\p3{}$'s.
\end{cor}

\begin{rem}\label{2.9}
  We already get a first interpretation in terms of abelian surfaces.
  Putting $W=W_3^\lor$ and choosing the unique
  $\p{}{}\sldfs$-invariant quartic $\K4$ of $\p{}{}W=\cp2+$ for $F$
  we get $V=W_7$, $\mathscr V_{2,9}'=\kappa (\p2+)$, $F^\flat=\Kp4$.
  Now if $a\in\p2+$ the proper transform of the quadric cone
  $\Gamma_a$ by $\kappa^{-1}$ is by \cite{GPo} birational to a Calabi
  Yau threefold, and by the preceding corollary contains --- when
  $\C{a}$ is smooth --- {\it two} dictinct pencils of special surfaces
  : the one induced by the six secant planes, parametrized by
  ${\mathcal C_a}$ and corresponding generically to abelian surfaces,
  and another one induced by the second ruling (parametrized by
  $\C{a}$) of planes of the cone. We think that these last ones are
  the same as the ones evoked in \cite[remark 5.7]{GPo}.
\end{rem}

\subsection{Boundary for the Klein quartic}

Note that the boundary for the general case is easily deducible from
what follows, but one need to introduce a covariant of $F$, this would
be beyond the subject of this paper.

So in this section, we focus on the surface $\Delta_F=\{\square \in
VSP (F,6)\, | \, \lambda_\square\neq (1^6) \}$ when the quartic $F$
admits $\p{}{}\sldfs$ as its group of automorphisms.

We start by choosing a faithful embedding $1\longrightarrow
\sldfs\longrightarrow \sltc$, so the vector space $W$ of our preceding
section becomes a $\sldfs$--module (necessarily irreducible), say
$W\simeq W_3^\lor$ and the decomposition $S^4W_3={\bf C}\oplus W_6
\oplus W_8$ allows us to consider the {\it unique}
$\p{}{}\sldfs$--invariant quartic $\K4$ of $\cp2+=\p{}{}W_3$. Such a
quartic is called a Klein quartic and becomes the quartic $F$ of our
preceding section.  All the quartic covariants of $F$ are equal to $F$
(when non zero) and the Klein quartic $\Kp4\subset\p2+$ is (by unicity)
the quartic $F^\flat$ of the last section.

We'll need the classical
\begin{lemma}
  There is a unique $\sldfs$--invariant even theta characteristic
  $\vartheta$ on the genus $3$ curve $\K4$ ({\it resp.} $\Kp4$).
\end{lemma}
\begin{preuve}
  The existence follows directly by the existence of a
  $\sldfs$--invariant injection $W_3\longrightarrow S^2U_4$ so that
  one can illustrate $\K4$ as the jacobian of a net of quadrics (in
  $\p{}{}U_4^\lor$). It is well known that such jacobian is endowed
  with an even theta characteristic (cf \cite{B}). Reciprocally, such
  a theta characteristic on a curve of genus $3$ comes with a net of
  quadrics and ${\rm Hom}_{\mathrm{SL} (2,{\bf F}_7)} (W_3,S^2 U)\neq
  0$ if and only if the four dimensional vector space $U$ equals
  $U_4^\lor$ as $\sldfs$--module.
\end{preuve}

We have
\begin{prop}\label{triangle}
  Let $a\in\Kp4$ and $(x_1,x_2,x_3)\in \Kp4\times\Kp4\times\Kp4$ such that
  ${\rm h}^0 (\vartheta +a-x_i)=1$, then the antipolar conic $\C{x_i}$
  of $x_i$ with respect to $\K4$ contains $x_j$.
\end{prop}
\begin{preuve}
  Let us leave the plane $\p2+$ and take a look at the configuration
  in $\p5+=\p{}{}S^2W_3=\p{}{}W_6=\cp5+$. The image of $x_i$ by the
  Veronese embedding $\nu_2$ lie on the quadric $Q_{\mathscr K_4'}$.
  On the other hand, noticing that ${\rm Hom}_{{\rm SL}_2{\bf F}_7}
  ({\bf C}, S^2S^2W_3)={\bf C}$ this quadric can be interpreted
  \begin{itemize}
  \item as the inverse of the quadric $Q_{\mathscr K_4}$;
  \item as the Plücker embedding of the Grassmannian of lines of
    $\p3-$ using the $\sldfs$-invariant identification $S^2W_3\simeq
    \Lambda^2U_4$.
\end{itemize}
Let us denote by $\Kp6$ the jacobian of the net of quadrics given by
$W_3\longrightarrow S^2U_4$ and remember that this curve is (by
unicity) canonically isomorphic to $\Kp4$ itself. So $\nu_2 (x_i)$ is a
line in $\p3-$ (still denoted by $\nu_2 (x_i)$) and this one turns out
to be a trisecant line to the sextic $\K6$ containing the image of $a$
by the identification $\Kp4=\Kp6$. This interpretation of the $(3,3)$
correspondence on $\Kp4=\Kp6$ induced by the even theta characteristic
as the incidence correspondence between $\Kp6$ and its trisecant lines
is due to Clebsch. Now the three lines $\nu_2 (x_i)$ are concurrent in
$a$ and then the three points $\nu_2 (x_i)$ of $\p5+$ span a
projective plane contained in the inverse of the quadric $Q_{\mathscr
  K_4}$.  In particular, $\alpha_{{\mathscr K}_4}^{-1}(\nu_2 (
x_i))\cdot \nu_2(x_j)=0$ which is precisely what we need to claim that
$x_j \in \C{x_i}$.
\end{preuve}

Notice that using the same geometric interpretation we get immediately
\begin{cor}
  If $a\in\Kp4$ then the antipolar conic $\C{a}$ intersects the hessian
  triangle $T_a$ ({\it ie} the hessian of the polar cubic of $a$ with
  respect to $\Kp4$) in points of the quartic $\Kp4$ (and $\C{a}\cap
  \Kp4-2a=T_a\cap \Kp4 -2 x_1-2x_2-2x_3)$).
\end{cor}

\begin{prop}\label{types}
  Let $p\in \Delta:=\Delta_{{\mathscr K}{\lower.5ex\hbox{\kern-.25em
        {\rm\tiny 4}}}}$, then there exists at least one point $a$ in
  the support of $\zeta_p$ such that $a\in\Kp4$ and the type of
  $\zeta_p$ is one of the following
  $$\hbox{\begin{tabular}{c|c|c}
      & $a\in\hspace{-8pt}/\hspace{4pt}\H6$ &\multicolumn{1}{c}{$a\in\H6$}\\
      \hline \multirow{3}{2cm}{type of $\zeta_p$}
      &$2,1,1,1,1$&$2,2,1,1$\\
      &$3,1,1,1$&$4,2$\\
      &$2,2,2$&$(2,2,2)_s$
\end{tabular}}$$
\end{prop}

\begin{preuvesc}
  From the preceding section, a point $p$ of $VSP(\K4,6)$ is in
  $\Delta$ if and only if the support of $\zeta_p$ intersects the
  quartic curve $\Kp4$.  So let $a\in\Kp4$, by corollary \ref{conics}
  the only thing to understand is the type of $\zeta_p$ when the point
  $p$ moves along the conic ${\mathcal C}_a$.  We have the
  alternative: the conic ${\mathscr C}_a$ is smooth (case {\it i}) or
  $a\in\H6:=H_{\mathscr K_4'}$ (case {\it ii}).
    \begin{enumerate}
    \item\label{types1} Denote (once again) by $S_n$ the two
      $(n+1)$-dimensional vector space ${\rm H}^0 (\O_{\mathcal C_{a}}
      (n))$, we have $S_n=S^nS_1$.  As $a\in \Kp4$, the $(1,5)$
      correspondence between the two (isomorphic) rational curves
      ${\mathcal C_a}$ and ${\mathscr C_a}$ has a base point, namely
      the point $a$ itself on ${\C{a}}$ and then reduces to a $(1,4)$
      correspondence.  The induced pencil of divisors of degree $4$ in
      $\p{}{}S_4$ intersects the variety of non reduced divisors in
      six points (as any generic pencil in $\p{}{}S_4$) and the
      expected types of $\zeta_p$ are hence $(2,1,1,1,1)$ generically,
      $(3,1,1,1)$ once and $(2,2,1,1)$ six times (each corresponding
      to a point of $\C{a}\cap\Kp4 -\{a\}$).  But by the preceding
      proposition, if $a'\in \C{a}\cap\Kp4$ and $a\neq a'$ then the two
      conics $\C{a}$ and $\C{a'}$ intersect in $a+a'+2a''$ with
      $a''\in\Kp4$ hence the six expected $(2,2,1,1)$ on $\mathcal C_a$
      become three $(2,2,2)$ for the particular Klein quartic. Notice
      that in such a case, the scheme $\zeta_p$ has a length
      decomposition $2\cdot(a+a'+a'')$ and there exists a point
      $b\in\Kp4$ so that ${\rm h}^0(\vartheta +b -x)=1$ whenever
      $x\in\{a,a',a''\}$. Let us denote by $b_x$ the intersection of
      $\C{x}$ with the line $\overline{xb}$, then $$a^3\cdot
      b_a+a'^3\cdot b_{a'}+a''^3\cdot b_{a''}=0 $$
      is an equation of
      $\K4$.
    \item suppose now the point $a$ is one of $24$ points of
      intersection of the quartic $\Kp4$ and its Hessian $\H6$. Such
      points come $3$ by $3$ and the group $\mu_3$ acts on each
      triplet (so there is an order $a_1,a_2,a_3$ on such triplet).
      Put $a=a_1$.  The conic $\C{a_1}$ is no longer smooth and
      decomposes in two lines, say $\lr=a_1a_2$ and
      $\lr'=a_2a_3$. 
      Each generic point $b$ of the line $\lr$ gives us a point
      $2a_1+b+b'+2a_3=\C{a_1}\cap\C{b}+a_1+b$ of $\Delta$ (hence of
      type $(2,2,1,1)$) with $b'\in\lr$ defined such that the degree
      $4$ divisor $a_2+a_1+b+b'$ on the line $\lr$ is harmonic.  One
      can even provide the corresponding equation of the quartic
      $\K4$

      \begin{center}
        $\displaystyle\epsilon(\beta x+\alpha z)^4-\epsilon(\beta
        x-\alpha z)^4- 2\alpha \beta\{(x+\epsilon(\beta^2z-\alpha
        ^2y))^4-x^4\}+2\alpha ^3\beta((y+\epsilon\, z)^4-y^4)=0$
    \end{center}
    
    with coefficients in ${\bf C}[\epsilon]/\epsilon^2$.  We can
    forget points of $\Delta$ arising from a point of $\lr'$, for such
    points can be constructed as the preceding ones by starting with
    the point $a_2$ instead of the point $a_1$. The possible
    degeneracies follow easily: when $(\alpha:\, \beta)$ tends to
    $(1:\, 0)$ we get back to the well known $(2,2,2)_s$ case, the
    last equation becomes
    $$(z+\epsilon\, x)^4-z^4+(x+\epsilon\,y)^4-x^4+(y+\epsilon\,
    z)^4-y^4=4\epsilon\, (z^3x+x^3y+y^3z) =0.$$
    The last possible
    degeneration arises when $(\alpha:\, \beta)$ tends to $(1:\, 0)$
    in which case we get $(4,2)$ as partition of $6$.\hfill $\square$
    \end{enumerate}
\end{preuvesc}


\section{Degenerated abelian surfaces}

Now that we understood what was the boundary $\Delta$ of $V_{22}$ in
${\rm Hilb}(6,\p2+)$ we try to understand it in terms of what we call
`degenerated abelian surfaces'. The method we are going to employ is
very na\"ive : given $s\in \Delta$, find a surface $A_s$ in $\p{}6$
which intersects $\p2+$ along $s$ and check that this surface is sent
to a projective plane by the map $\kappa$, {\it ie} that ${\rm
  h}^0(\I_{A_s} (7))^{G_7}=5$. Obviously translation scrolls are
candidates and we'll see these are the good ones.  {In this section,
  we work up to the action of $\p{}{}\sldfs$.}

\subsection{Translation scrolls}
We'll need the
\begin{prop}\label{2tscroll}
  Every translation scroll of an elliptic normal curve of degree $7$
  by a 2-torsion point is a smooth elliptic scroll of degree $7$ and
  contains 3 elliptic normal curves of degree $7$.
\end{prop}
\begin{preuve}
  cf. \cite[proposition 1.1]{CH} or \cite{Mie}.
\end{preuve}





Let us start with $s\in \Delta_\lambda$ with
$\lambda\in\{(2,1,1,1),(3,1,1,1)\}$. Only one point of $s$ has a
multiple structure, say $p\in\Kp4$ and the support of $s$ consists in
four {\it distinct} points on the conic $\C{p}$. The stabilizer of
$(s)_{\rm red}$ under $\p{}{}\sldc\simeq {\rm Aut} (\C{p})$ is in
general isomorphic to ${\bf Z}_2^2$ (if it is bigger consider the
subgroup $\{{\rm Id},(1,2) (3,4),(1,3) (2,4), (1,4) (2,3)\}$ in ${\rm
  Stab}_{{\bf P}{\rm SL} (2,{\bf C})} ((s)_{\rm red})\subset {\mathfrak
  S}_4$).  Consider the double cover $E_s$ of $\C{p}$ ramified at the
four points ${\bf Z}_2^2\cdot p$. It is a smooth elliptic curve and we
choose $p$ as origin on $E_s$. Note that, by proposition \ref{triangle},
$E_s$ depends only of $p$ so it will be denoted by $E_p$.  Now the
linear system $|7\cdot p|$ is a $H_7$-module and we can embed $E_p$ in
$\p{}6$ in $168$ distinct ways, one of them send $p\in E_p$ to
$p\in\p2+$.  Next $p\in s$ is doubled along a line that will intersect
$\C{p}$ in a further point, say $p_s$ (of course $p=p_s$ if
$\lambda=(3,1,1,1)$).  Denote by $\sigma_s$ one of the inverse images
of $p_s$ by the $2:1$ map $E_p\longrightarrow \C{p}$.

\begin{prop}
  The translation scroll $(E_p,\sigma_s)$ intersects $\p2+$ along $s$
  and is mapped by $\kappa$ to a plane.
\end{prop}

\begin{preuve}
  First the bisecant variety of $E_p$ intersects $\p2+$ along $\C{p}$
  : indeed such variety intersects $\p2+$ along a conic ($E_p$ being
  of degree $7$ and invariant under a symmetry which preserves
  $\p2+$). Next by the previous proposition such conic must contains
  the three pairs of points of $\Kp4^2$ such as $(q_1,q_2)$ where
  $\{p,q_1,q_2\}$ are associated to a same point of $\Kp4$ under the
  $\p{}{}\sldfs$-invariant $(3,3)$ correspondence on $\Kp4$. By
  proposition \ref{triangle}, this conic is nothing but $\C{p}$.
  
  Next as $\pm\sigma$ moves along $\C{p}$, the set $(E_p,\pm
  \sigma)\cap \C{p}-{2p}$ describes a {\it pencil} of degree $4$
  divisors on $\C{p}$, we need to identify this pencil with our pencil
  of degree $4$ divisors on $\C{p}$ (proof of proposition \ref{types},
  item \ref{types1}). But both pencils contain the three divisors of
  type $(2,2)$ such as $2q_1+2q_2$ so they are equal.
  
  The point now is to prove that $(E_p, \sigma_s)$ is mapped to a
  plane under $\kappa$ {\it i.e.} that ${\rm h}^0 (\O_{(E_p,
    \sigma_s)} (7))^{G_7}=3$ or rather ${\rm h}^0 (\O_{(E_p,
    \sigma_s)} (7))^{G_7}\leqslant 3$ by the previous paragraphs. The
  scroll $(E_p, \sigma_s)$ is a $\p{}1$--bundle over $E_p$ (in two
  ways, these correspond to the choices $\sigma_s$ and $-\sigma_s$ to
  define the scroll) so we have a map $$(E_p, \sigma_s)\longrightarrow
  E_p$$
  and if $R$ denotes a generic fiber we get a sequence
  $${\rm H}^0 (\O_{(E_p, \sigma_s)}(7))^{G_7}\longrightarrow {\rm H}^0
  (\O_R(7))^{G_7}\longrightarrow 0$$
  which turns out to be exact: if a
  $G_7$--invariant \septic{} $S$ contains the line $R$, then its
  intersection with $(E_p, \sigma_s)$ contains $E_p\cup G_7\cdot R$
  which is of degree $7+49\times 2\times 1$. Now our scroll is of
  degree $14$ and by Bezout we conclude $(E_p, \sigma_s)\subset S$.
  Using a semi-continuity argument we just need to find one fiber $R$
  such that ${\rm h}^0 (\O_R(7))^{G_7}\leq 3$. Choose one of the
  forty-nine $\p2+$ and pick up one of the four lines of $(E_p,
  \sigma_s)$ which intersects it. Then the restriction
  $${\rm H}^0 (\O_{{\bf P}V_0}(7))^{G_7}\longrightarrow {\rm H}^0
  (\O_R(7))^{G_7}$$
  is of rank $3$ at most which is precisely what we
  needed.
\end{preuve}

Of course to get the $(2,2,2)$ cases it is natural to make $\sigma_s$
tend to a $2$-torsion point. Then the translation scroll $(E_p,
\sigma_s)$ tends to a {\it smooth} scroll of degree $7$ and everything
is lost in virtue of the

\begin{rem}
  Any translation scroll $(E_p, \sigma_s)$ where $\sigma_s$ is a
  non-trivial $2$-torsion point of $E_p$ is contained in all our
  $G_7$-invariant \septic{} hypersurfaces. Indeed such a scroll
  intersects any of the forty-nine $\p3-$'s along a line and contains
  the curve $E_p$ itself, once again Bezout together with the
  inequality $49\times 1+7>7\times 7$ allow us to conclude. However we
  have
\end{rem}

\begin{prop}
  If $\lambda_s\in (2,2,2)$ then there exist an elliptic curve $E_p$,
  a two torsion point $\sigma_s$ on $E_p$ and a double structure on
  the translation scroll $(E_p, \sigma_s)$ intersecting $\p2+$ along
  $s$ and mapped to a plane by $\kappa$.
\end{prop}

\begin{preuve}
  Let $s\in \Delta_{(2,2,2)}$. By proposition \ref{2tscroll} one and
  only one smooth translation scroll $X$ intersects $\p2+$ along
  $(s)_{\rm red}$ so we just need to find a double structure
  $\tilde{X}$ on $X$ such that ${\rm h^0} (\I_{\tilde{X}}
  (7))^{G_7}=5$.  Now $X$ contains two (in fact three by
  \ref{2tscroll}) elliptic curves $E_p$ and $E_{p'}$ and by definition
  is contained in the two corresponding bisecant varieties $S_{E_{p}}$
  and $S_{E_{p'}}$. But these two varieties are the proper transforms
  by $\kappa^{-1}$ of the two quadric cones $\Gamma_p$ and
  $\Gamma_{p'}$ (cf corollary \ref{CY}). These two cones intersect
  along the six secant plane to $\kappa (\p2+)$ corresponding to $s$
  so we are done!  The double structure is then easy to understand:
  one considers the double structures on $X\backslash E_p$ ({\it
    resp.}  $X\backslash E_{p'}$) defined by the embedding
  $X\longrightarrow S_{E_{p}}$ ({\it resp.}  $X\longrightarrow
  S_{E_{p'}} $) and such structures coincide on $X\backslash (E_p\cup
  E_{p'})$.
\end{preuve}

\subsection{Union of seven quadrics}\label{the7p3s}

We still have to consider the missing cases, namely schemes $s$ such
that $\lambda_s\in\{(2,2,1,1),(2,2,2)_s,(4,2)\}$.  These are
degenerations of the preceding ones. A degenerated elliptic curve is
nothing but a heptagon, and such curves come by triplets
$(E_0,E_1,E_2)$ (fig. 3) with $E_i= \bigcup_{k=0}^6 \overline{e_k
  e_{k+1+i}}$ and $\{e_x\}_{x\in {\bf Z}_7}$ is an orbit of minimal
cardinality under the action of $G_7$.

The Heisenberg action on each curve $E_i$ reduces to an action of
${\bf Z}_7$. Let us denote by $\p{I}{}$ the projective space spanned
by the points $\{e_i\}_{i\in I}$.

For $i\in \{0,1,2\}$ put $B_i=\bigcup_{k=0}^6 \p{I_i^k}{}$ with
$I_i^k=\{k+i+1,k-i-1,k+3i+3,k-3i-3\}$. Then the bisecant variety of
$E_i$ is $B_j+B_k$ with $\{i,j,k\}=\{0,1,2\}$.

Finally let us choose one of the forty nine \p2+ and suppose $E_i$
intersects it in $a_{i+1}$. We are then ready to check the remaining
cases:

Let $s\in \Delta$ such that $\lambda_s= (2,2,1,1)$ and $s=2\cdot
a_1+2\cdot a_3+b+b'$ with $b$ and $b'$ on the line $\overline{a_1a_2}$
(Fig. 1). So we have $s\in {\mathcal C}_{a_1}\cap {\mathcal C}_{a_3}$.
The corresponding degenerated abelian surface $A_s$ needs to be on the
bisecant varieties of $E_0$ and $E_2$ so we have $A_s\subset B_1$. As
$B_1$ is the union of seven $\p{}{3}$'s the surface $A_s$ is the union
of seven quadrics.

When $s$ moves along ${\mathcal C}_{a_1}\cap {\mathcal C}_{a_3}$ we
get the two other kinds of degenerations :

\begin{itemize}
\item if $\lambda_s= (2,2,2)_s$, then the surface $A_s$ degenerates in
  the union of the fourteen planes $B_0\cap B_1\cup B_1\cap B_2\cup
  B_2\cap B_0$;
\item if $\lambda_s= (4,2)$, then the surface $A_s$ degenerates in the
  union of the seven planes $\bigcup_{k=0}^6\p{\{k,k+3,k-3\}}{}$
  doubled along $B_1$.
\end{itemize}

\subsection{The smooth case}\label{ANTI} 
To begin with let us

\begin{rem}\label{3.5}
  Choose coordinates $(y_i)_{i\in \{0,...,6\}}$ in $V_0$ together with
  one of the forty-nine $\p3-$'s of equations $y_4=y_3,\ y_5=y_2,\ 
  y_6=y_1$.
  Using the $N_7$-invariant isomorphism $\Lambda^2V_4=W_3^\lor\otimes
  V_0$ let us introduce for $x= (x_1:x_2:x_3)\in \p{}{}W_3$ and $y=
  (y_0:-:y_4)\in \p3-\subset\p{}{}V_0$ the matrix
  $$M_y(x)= \left(\matrix
    {x_2y_2&-x_3y_1-x_1y_3&x_3y_0&-x_2y_1-x_1y_2\cr
      -x_3y_3&-x_3y_2+x_2y_3&-x_2y_1+x_1y_2&x_1y_0\cr
      -x_1y_1&x_2y_0&x_3y_1&x_3y_2+x_2y_3}\right).
  $$
  It is the restriction to $\p3-$ of a skew-symmetric Moore matrix
  $M (x,y)$ (see remarks \ref{smoore}, \ref{2.9} and \cite{GPo}). This
  matrix $M (x,y)$ defines the Calabi Yau threefold which is the
  strict transform of the cone $\Gamma_x$ by $\kappa^{-1}$.  For
  general $y\in \p3-$, $M_y (x)$ defines $6$ points in $\p{}{}W_3$,
  meaning there is {\it one} abelian surface containing $y$ and six
  Calabi Yau
  threefolds of the preceding type which contain this surface. But $M_y (x)$ may degenerate for special points $y\in \p3-$ (for such cases we get more than six points in $\p{}{}W_3$):\\
  \begin{center}
\begin{tabular}{|c|c|c|c|}
   \hline 
   rank of $\Lambda^3M_y (x)$& $y$ in & $x$ in & abelian varieties passing through $y$\\
   \hline 3 & ${\mathscr K}_6$ & ${\mathscr H}_6$ & have a 3--secant line\\
   \hline 2 & $C_{18}$ & ${\mathscr K}'_4$ & translation scrolls\\
   \hline 1 & $Z$ &${\mathscr K}'_4\cap {\mathscr H}_6$ & reducible\\
   \hline
  \end{tabular}
\end{center}
where ${\mathscr K}_6$ is the unique $\p{}{}\sldfs$-invariant curve of
degree $6$ and genus $3$ in $\p3-$, $C_{18}$ is a
$\p{}{}\sldfs$-invariant curve of degree $18$ and genus $35$ in $\p3-$
(analogue of the Bring curve in the $(1,5)$ case) and $Z$ is the
minimal orbit (of cardinality eight) under the action of
$\p{}{}\sldfs$ on $\p3-$.
\end{rem}

So far, and summing up all the results of this section, we need, in
order to complete the proof of the theorem, to show that $A$ is smooth
and abelian provided the type of $\zeta_A$ is $(1,1,1,1,1,1)$.  Now by
the Enriques-Kodaira classification of surfaces (see \cite[chapter
VI]{bpvdv}) complex tori are entirely characterized by their numerical
invariants. So any generalized $G_7$--embedded abelian surface is an
abelian surface provided it is smooth. We begin by the

\begin{lemma}
  A generalized $G_{7}$-embedded abelian surface $A$, singular along a
  curve, intersects `the' plane $\p2+$ in a non reduced scheme.
\end{lemma}
\begin{preuve}
  Let $A$ a singular generalized $G_{7}$-embedded abelian surface. The
  proposition is obviously true if $A$ is singular in codimension $0$
  that is to say if $A$ carries a double structure. For such surface,
  the intersection of its reduced structure (of degree $7$) with any
  $\p2+$ cannot be six distinct points so $\zeta_A=A\cap \p2+$, which
  is of length six, cannot be reduced.
  
  By assumption the singular locus of $A$ contains a curve $C$.  We
  can also assume $C$ is $G_7$-invariant (if not we replace $C$ by its
  orbit under $G_7$).

  
  If $C$ has degree $7$, it is necessarily elliptic and, being $G_7$
  embedded, intersects $\p2+$ in a point of the Klein curve $\Kp4$ so
  we are done. Indeed the rationality of $C$ (if irreducible) is
  totally excluded (such curve admits either a unique 4-secant plane,
  a unique trisecant line or a (unique) double point, this would be a
  contradiction with the irreducibility of $V_0$ as $H_7$--module),
  but $C$ can still splits in the union of seven lines. We want to
  prove that $C$ is a heptagon (that is to say elliptic).  The
  stabilizer of one of the lines under the action of $H_7$ is
  isomorphic to ${\bf Z}_7$ so we get, on each line $\lr \subset C$,
  two fixed points under the action of ${\rm Stab}_{H_7} (\lr)\simeq
  {\bf Z}_7$ and then an orbit of fourteen points on $C$. Noticing all
  the components have the same stabilizer (the only group morphism from
  ${\bf Z}_7$ to the symmetric group ${\mathfrak S}_6$ is constant)
  and considering the symmetries of $G_7$ it is easy to prove that
  these fourteen points coincide two by two, implying $C$ is a
  heptagon.
  
  If the singular locus $C$ has degree $14$, then the reduced
  structure of it has degree $7$ or $14$.  Only the latter is a
  problem.  The normalization of the surface has sectional genus
  $(-6)$ so it consists of at least seven components. They have the
  same degree {\it i.e.} $1$ or $2$, so their number must be $7$ and
  their degree must be $2$. In particular, either the surface $A$ is
  contained in one orbit of seven $\p3{}$'s under $G_7$ or the reduced
  structure of $A$ consists of seven planes. In both cases, it is a
  simple matter to conclude for the only orbits of seven $\p3{}$'s
  under $G_7$ are listed in the subsection \ref{the7p3s} (consider for
  instance their possible instersections with the forty-nine $\p3+$)
  so the surface $A$ appears already in the subsection \ref{the7p3s}
  and the proposition is true for such surfaces.
  
  The last possible case is when $C$ has degree $21$, but then the
  surface $A$ has $14$ components (its normalization has sectional
  genus $(-13)$) so $C$ splits and, as $\gcd (49,21)=7$, contains
  three $G_7$-invariant curves of degree $7$ so we are back to the
  first case.
\end{preuve}

\noindent{\it End of the proof of the theorem.}


\noindent Let $A$  be a generalized abelian surface, preimage of a six
secant plane of $\kappa (\p2+)$ by $\kappa$, \textit{i.e.}
$\zeta_A=(1,1,1,1,1,1)$.  Since $A\cap \Kp4 =\emptyset$ we may divide
in two cases; $A$ intersects ${\mathscr H}_6$ but not $\Kp4$ in
$\p2+$, and $A$ intersects neither ${\mathscr H}_6$ not $\Kp4$.

\begin{enumerate}
  
\item Assume $A$ intersects ${\mathscr H}_6$ in $\p2+$, then $A$ is a
  {\it smooth} plane curve fibration and has a trisecant line in
  $\p2+$: see construction in \cite{HR};
\item Assume $A\cap {\mathscr H}_6=A\cap {\mathscr K}'_4=\emptyset$.
  Otherwise $A$ is a translation scroll or a plane curve fibration.
  By the previous lemma, $A$ is irreducible with isolated
  singularities. Let us choose $a\in \p2+\cap A$ and consider
  (identifying $\p2+$ with $\p{}{}W_3$) the Calabi Yau threefold $Y_a$
  preimage of the cone $\Gamma_a$ by $\kappa$. We know that
  $Y_a=\bigcup_{t\in {\mathcal C}_a} A_t$. $Y_a$ has a quadratic
  singularity at $a$, and the general surface $A_t$ is smooth at $a$,
  so after a small resolution $A_t$ will be Cartier there.  In the
  following we will ignore this distinction between $Y_a$ and its
  small resolution of $a$ and its $G_7$ translates.  We have
  \begin{itemize}
  \item $\omega_A={\mathscr O}_A$. Indeed, the general $A_t$ is smooth
    (since it is for $t\in {\mathscr H}_6$, notice also that
    $\bigcup_{t\in {\mathcal C}_a} \{A_t \} $ contains $8$ translation
    scrolls and $4$ elliptic curve fibrations) with normal bundle
    ${\mathscr O}_A (A)={\mathscr O}_A$ so $\omega_A={\mathscr O}_A
    (K_{Y_a})={\mathscr O}_A$;
  \item $A$ is a Cartier divisor on the Calabi Yau threefold $Y_a$:
    the threefold $Y_a$ is mapped by definition to a quadric
    $\Gamma_a$ of rank $4$ (this follows easily considering $a$ is not
    on ${\mathscr H}_6$) so $A$ is a Cartier divisor on $Y_A$
    providing $\bigcap_{t\in {\mathcal C}_a}=G_7\cdot a$ {\it i.e.} if
    the cardinality of $\bigcap_{t\in {\mathcal C}_a}$ is $49$.
    
    Let $B_a=\bigcap_{t\in {\mathcal C}_a}$. If $B_a$ contains a curve
    $D$ then it has degree at least $14$.  On the other hand it must
    lie on a the translation scrolls of the pencil $A_t$. Let $T$ be
    such a scroll and $\widetilde{T}$ its desingularization. Then
    ${\rm Pic} (\widetilde{T})=<E_0, F>$, where $E_0$ is a section
    with $E_0^2=0$ and F is a fibre. The preimage of the singular
    curve on $T$ is the union of two sections $E_0$ and $E_0'$.  Now,
    $B_a$ cannot intersect the singular curve on $T$, for degree
    reasons: $A_t$ is cut out by cubics, while the intersection would
    have cardinality at least $49$.  Therefore $D\cdot E_0=0$.  But
    then $D=\alpha\cdot E_0 +\beta\cdot E_0'$, which is again a
    contradiction.  So $B_a$ is finite.
    
    First $B_a\cap \p3- = \emptyset$: by remark \ref{3.5}, the points
    of $\p3-$ contained in a pencil of abelian surfaces are perfectly
    idenfified --- as well as the corresponding abelian surfaces:
    these are either translation scrolls or plane curve fibrations
    ---. Assume $B_a$ has cardinality greater than $49$, then because
    of the previous argument it has cardinality greater than
    $147=49+2\times 49$ where $2\times 49$ is the cardinality of $G_7$
    (acting on $\p{}{}V_0$).  But $B_a$ is contained in the base locus
    $B$ of the map $\kappa$ so $\kappa |_{A\ smooth}$ has at least
    $147+10\times 49$ ($147$ contained in $B_a$ and the other ones
    coming from the intersection of $A$ with the $49$ \p3-'s) base
    points, another contradiction.
  \end{itemize}
\end{enumerate}
\hfill{$\square$}

\section{Appendix}
\subsection{Some pictures}
\input{fig.degen} \input{fig.strat} \input{fig.lamelie}

\subsection{Some questions}
Let $A$ be a smooth abelian surface embedded in $\p{}{}V_0$. As
already noticed in proposition \ref{suite} we have a factorization
$$\xymatrix{{A^\flat}\ar[r]^{49:1\ }&A^{\lor\flat} \ar[r]^{2:1\ }
  &K_{A^\lor}^\flat\ar[r]^{2:1\ }&\kappa (A)}.$$
{What is the ramification locus of the last map?} Of course, the
map to the right needs to define a $K3$ surface so it is certainly a sextic
curve (of genus $4$), but it doesn't explain how to recover it (using
representation theory for instance).  A good understanding of this
point should enable one to get a reconstruction method, as in the
$(1,5)$ case, of the abelian surface $A^\lor$. Note that this
ramification locus admits the six points $\kappa (A)\cap \kappa
(\p2+)$ as double points.  Next the surface $A$ intersects $\p3-$ in
ten points, so projecting $A$ from $\p3-$ we get maps
$$\xymatrix{A^{\flat} \ar[r]^{2:1\ } &K_{A}^\flat\ar[r]^{2:1\ 
  }&\p2+}.$$
The last map is ramified along a sextic curve with the
six points $A\cap \p2+$ as double points.  Hence we have two maps
$\xymatrix{K_{A}^\flat\ar[r]^{2:1\ }&\p2+}$ and
$\xymatrix{K_{A}^\flat\ar[r]^{2:1\ }&\kappa (A^\lor)}$.

{\it Is it possible to find an (in fact $168$) indentification(s)
  between $\p2+$ and $\kappa (A^\lor)$ such that the two maps
  coincide?}  The answer is positive in the $(1,5)$ case and allows
one to identify the moduli space of $(1,5)$ polarized abelian surfaces
(without level structure) up to duality.  Note that one can easily
show that the two sets of six points $A\cap \p2+$ and $\kappa (A)\cap
\kappa (\p2+)$ are associated in Coble's sense. This phenomenom is in
fact true for any Fano threefold $V_{22}$: given a six secant plane to
${\mathscr V}_{2,9}'$ (Veronese surface isomorphic to $\p{}{}W$), it
intersects it along six points associated to the corresponding set in
$\p{}{}W$.
  
It is possible to show that the Fano threefold $VSP(\K4,6)$ is
`stable' by association of points {\it i.e.} that there exists a
(dual) Klein curve $C$ in the plane $\kappa (A)$ such that $\kappa
(A)\cap \kappa (\p2+)$ is a point of the corresponding variety $VSP
(C,6)$, isomorphic, up to $\p{}{}\sldfs$, to $VSP(\K4,6)$.  Is the
quotient the moduli space of $(1,7)$ polarized
abelian surfaces (without level structure) up to duality?



\end{document}

%% file: fig.degen.tex
\begin{minipage}{5cm}
\begin{center}
\psfrag{A}{$a$}
\psfrag{Ca}{$\C{a}$}
\includegraphics*[height=3.5cm ]{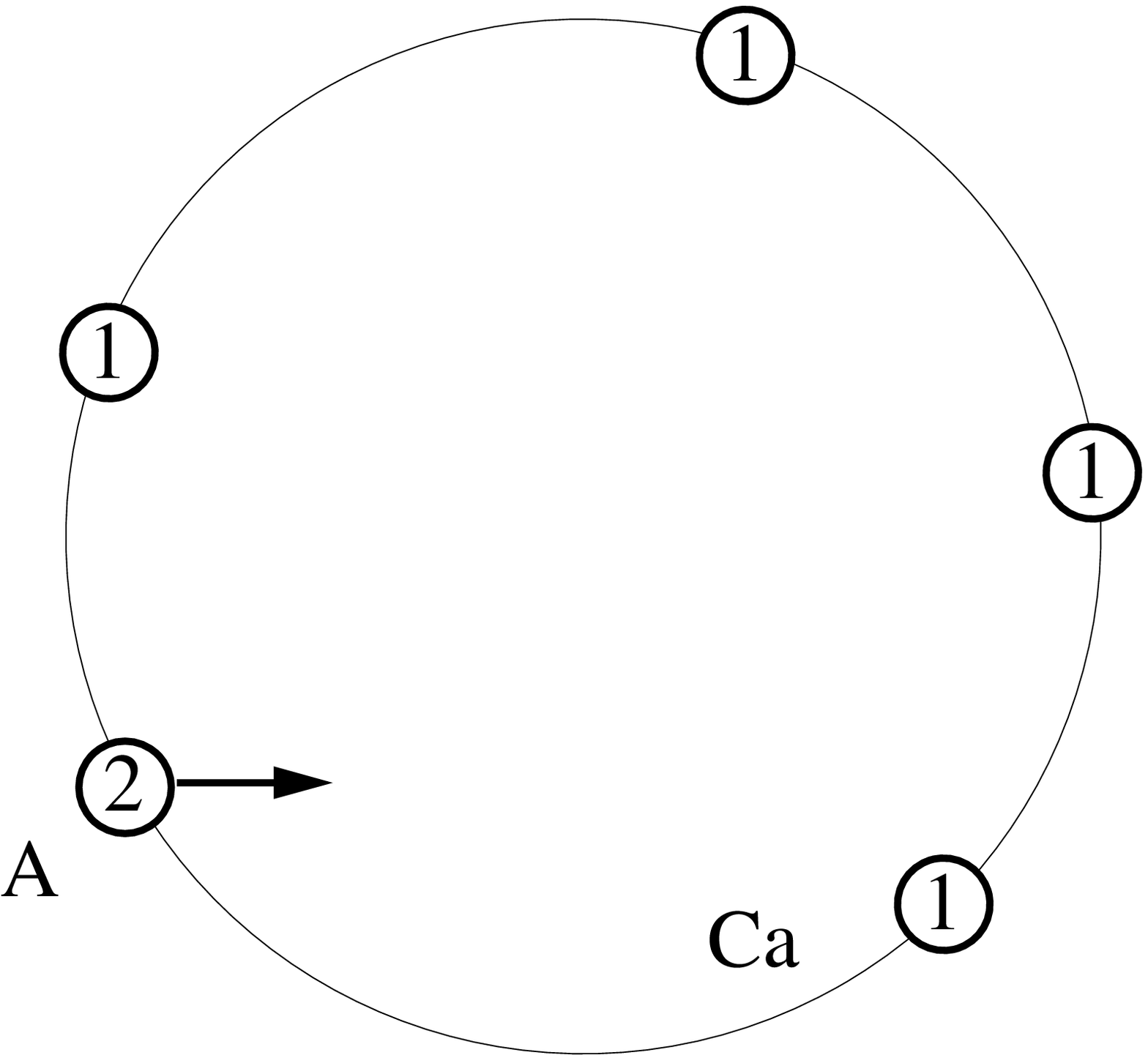}
  \label{fig:21111}
  \center{$\lambda_s=(2,1,1,1,1)$}
\end{center}
\end{minipage}\hfill
\begin{minipage}{5cm}
\begin{center}
\psfrag{A}{$a$}
\psfrag{Ca}{$\C{a}$}
\includegraphics*[height=3.5cm ]{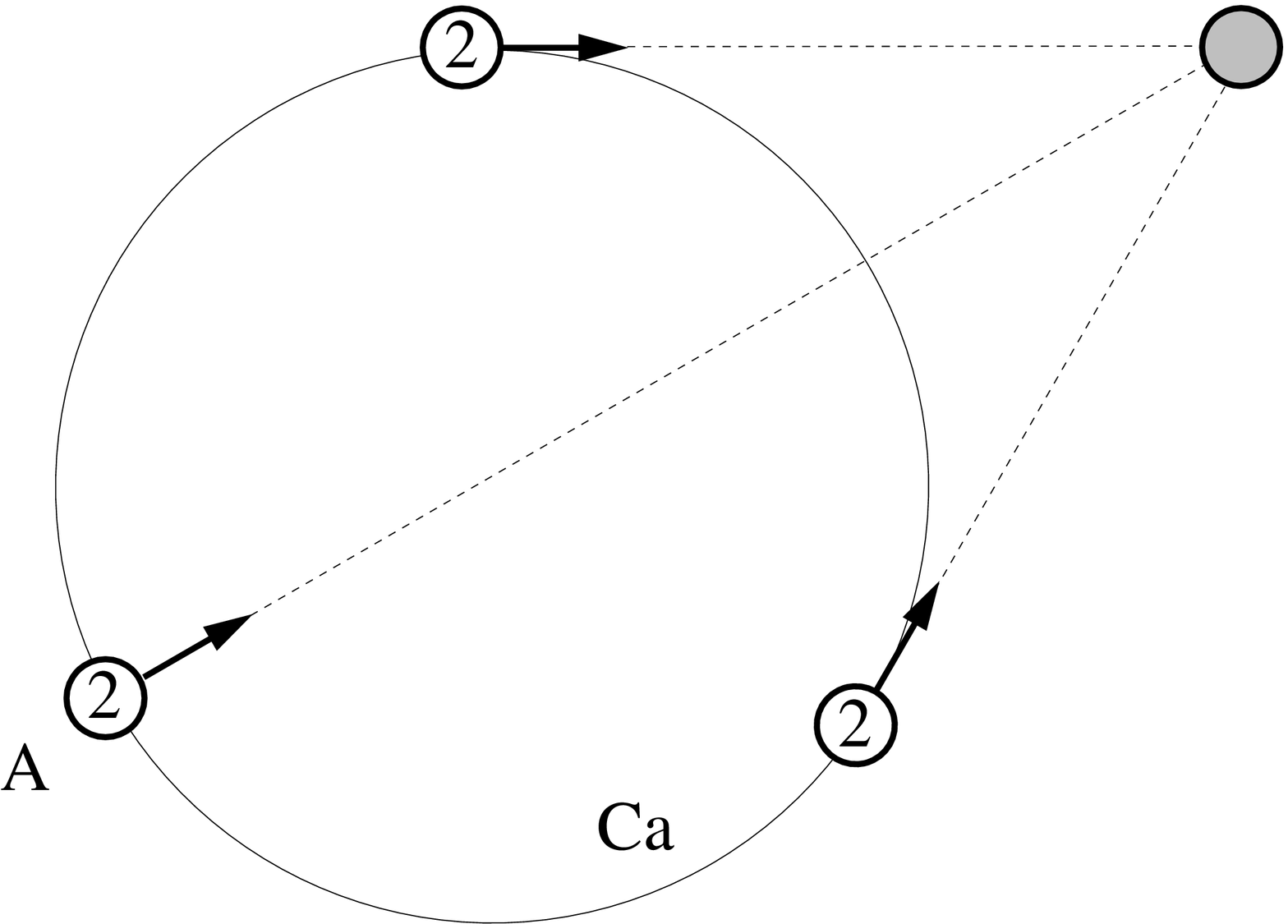}
  \label{fig:222gen}
  \center{$\lambda_s=(2,2,2)$}
\end{center}
\end{minipage}\hfill
\begin{minipage}{5cm}
\begin{center}
\psfrag{A}{$a$}
\psfrag{Ca}{$\C{a}$}
\includegraphics*[height=3.5cm ]{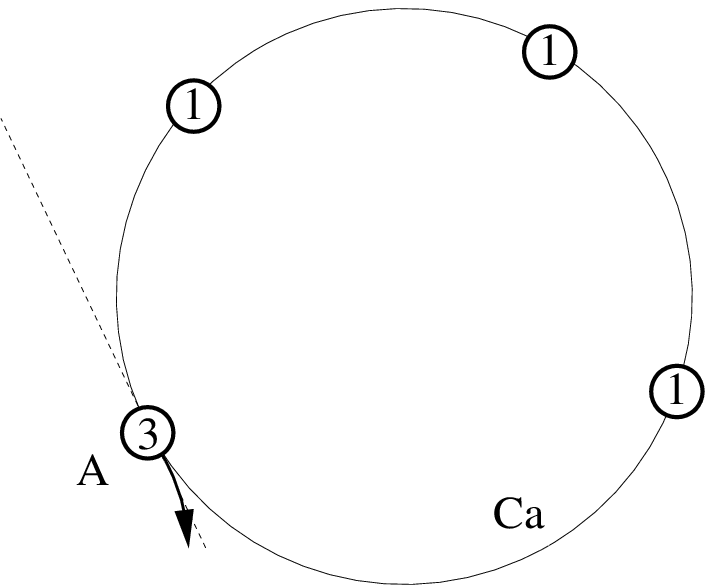}
  \label{fig:3111}
  \center{$\lambda_s=(3,1,1,1)$}
\end{center}
\end{minipage}
\vspace{1cm}

\hfill\begin{minipage}{5cm}
\begin{center}
  \psfrag{a1}{$a_1$}
  \psfrag{a2}{$a_2$}
  \psfrag{a3}{$a_3$}
  \includegraphics*[height=3.5cm ]{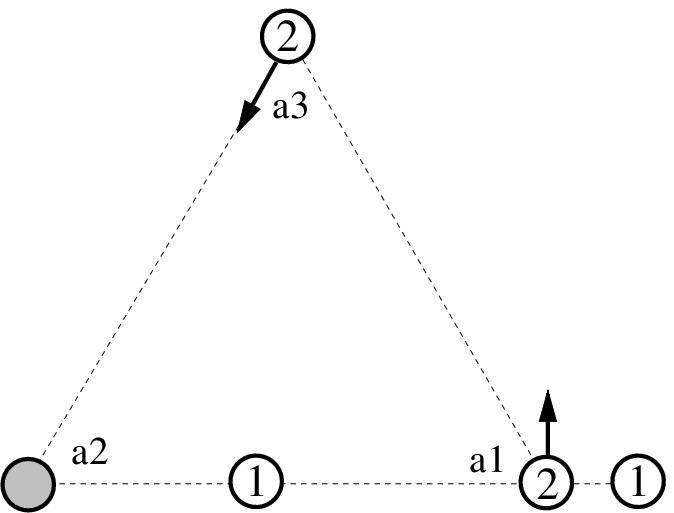}
  \center{$\lambda_s=(2,2,1,1)$}
  \label{fig:2211gen}
\end{center}
\end{minipage}\hfill
\begin{minipage}{5cm}
\begin{center}
  \psfrag{a1}{$a_1$}
  \psfrag{a2}{$a_2$}
  \psfrag{a3}{$a_3$}
  \includegraphics*[height=3.6cm ]{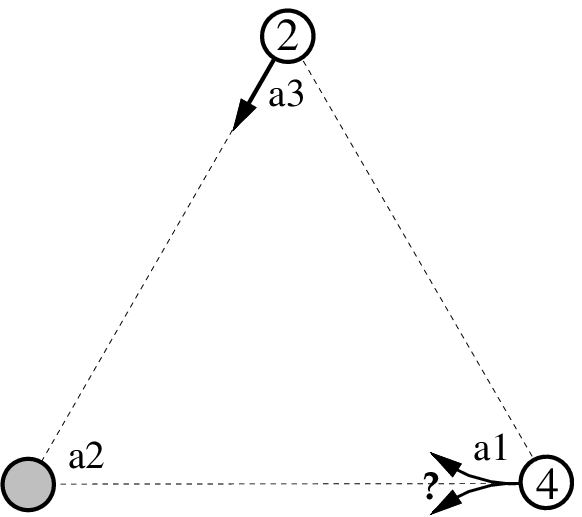}
  \center{$\lambda_s=(4,2)$}
  \label{fig:42}
\end{center}
\end{minipage}\hfill
\begin{minipage}{5cm}
\begin{center}
  \psfrag{a1}{$a_1$}
  \psfrag{a2}{$a_2$}
  \psfrag{a3}{$a_3$}
  \includegraphics*[height=3.5cm ]{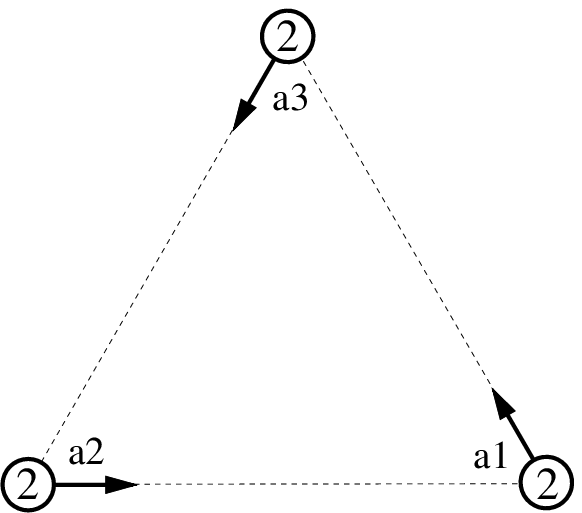}
  \center{$\lambda_s=(2,2,2)_s$}
  \label{fig:222x8}
\end{center}
\end{minipage}\hfill

\caption{
\label{degen} Possible configurations of $\zeta_s$ when $s\in \Delta$, where each arrow 
gives the (first) direction along which the point is doubled (oriented in a
purely decorative way).
}



%% file: fig.strat.tex
\begin{center}
  \psfrag{K2}{$\K4(2,2,2)$}
  \psfrag{K3}{$\K4(3,1,1,1)$}
  \psfrag{42}{$(4,2)$}
  \psfrag{222}{$(2,2,2)$}
  \psfrag{222S}{$(2,2,2)_s$}
  \psfrag{3111}{$(3,1,1,1)$}  
  \psfrag{2211}{$(2,2,1,1)$}
  \psfrag{21111}{$(2,1,1,1,1)$}
  \psfrag{b}{$b$}
  \psfrag{Cb}{${\mathcal C_b}$}
  \includegraphics*[height=7.5cm ]{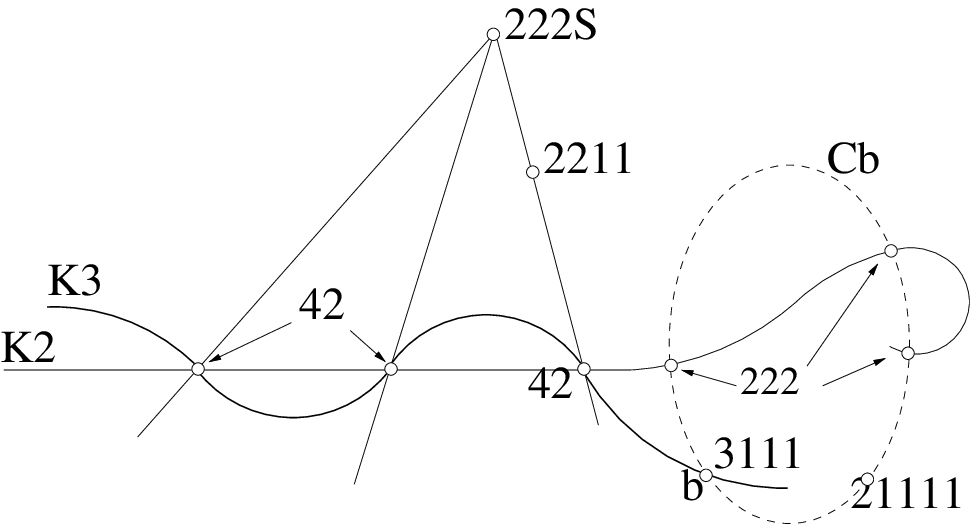}
\end{center}
\vspace{-1cm}\hspace{3.2cm}$\underbrace{\hbox{\hspace{7.1cm}}}_{\hbox{$8$ times}}$
\vspace{.5cm}
\caption{Stratification of the surface $\Delta$.}


%% file: fig.lamelie.tex
\begin{center}
  \psfrag{E0}{$E_0$} 
  \psfrag{E2}{$E_1$} 
  \psfrag{E1}{$E_2$} 
  \psfrag{e0}{$e_0$} 
  \psfrag{e1}{$e_1$} 
  \psfrag{e2}{$e_2$} 
  \psfrag{e3}{$e_3$} 
  \psfrag{e4}{$e_4$} 
  \psfrag{e5}{$e_5$} 
  \psfrag{e6}{$e_6$} 
  \psfrag{a1}{$a_1$} 
  \psfrag{a3}{$a_2$} 
  \psfrag{a2}{$a_3$} 
  \includegraphics*[height=6cm ]{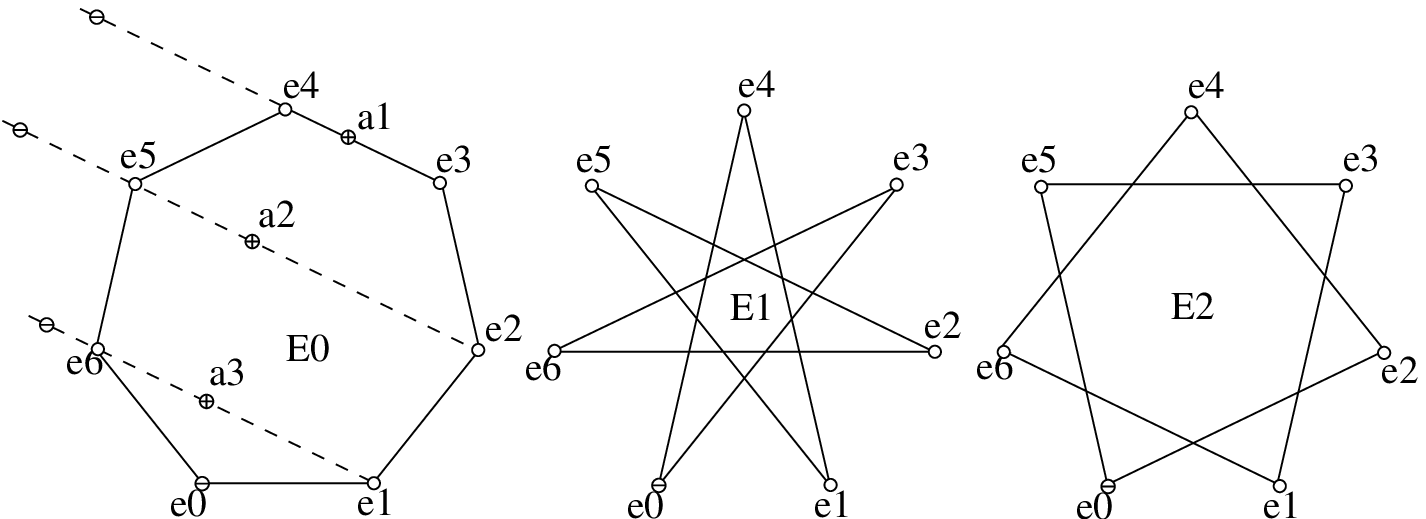}
  \caption{A triplet of degenerate elliptic curves.}
  \label{deluxe}
\end{center}


%% file: deg17.092003.bbl
\begin{thebibliography}{10}

\bibitem{B}
W.~Barth.
\newblock Moduli of vector bundles on the projective plane.
\newblock {\em Invent. Math.}, 42:63--91, 1977.

\bibitem{BHM}
W.~Barth, K.~Hulek, and R.~Moore.
\newblock Degenerations of the {H}orrocks-{M}umford surfaces.
\newblock {\em Math. Ann}, 277:735--755, 1987.

\bibitem{bpvdv}
W.~Barth, C.~{P}eters, and A.~Van de~{V}en.
\newblock {\em Compact complex surfaces}.
\newblock Number~4 in Ergebnisse der {M}ath. Springer-Verlag, 1984.

\bibitem{CH}
C.~Ciliberto and K.~Hulek.
\newblock A remark on the geometry of elliptic scrolls and bielliptic surfaces.
\newblock {\em Manuscripta Math.}, pages 231--224, 1998.

\bibitem{Dol}
I.~Dolgachev.
\newblock Invariant stable bundles over modular curves ${X}(p)$.
\newblock In {\em Recent progress in algebra (Taejon/Seoul, 1997)}, volume 224
  of {\em Contemp. Math.}, pages 65--99. Amer. Math. Soc., 1999.

\bibitem{D-K}
I.~Dolgachev and V.~Kanev.
\newblock Polar covariants of cubics and quartics.
\newblock {\em Advances in Math.}, 98:216--301, 1993.

\bibitem{Ell}
E.B. Elliot.
\newblock {\em An introduction to the Algebra of Quantics}.
\newblock Chelsea reprint, 1964.

\bibitem{Fr}
R.~Fricke.
\newblock {\em Lehrbuch der {A}lgebra {II}}.
\newblock Braunschweig, 1926.

\bibitem{F-H}
W.~Fulton and J.~Harris.
\newblock {\em Representation Theory}.
\newblock GTM 129. Springer-Verlag, 1991.

\bibitem{GPo}
M.~Gross and S.~Popescu.
\newblock Calabi-{Y}au threefolds and moduli of abelian surfaces {I}.
\newblock {\em Compositio Math.}, 127:169--228, 2001.

\bibitem{HR}
K.~Hulek and K.~Ranestad.
\newblock Abelian surfaces with two plane cubic fibrations and calabi-yau
  threefolds.
\newblock In Th. Peternell and F.-O. Schreyer, editors, {\em Complex Analysis
  and Algebraic Geometry}, pages 275--316. de Gruyter, 2000.

\bibitem{LB}
H.~Lange and C.~Birkenhake.
\newblock {\em Complex Abelian Varieties}.
\newblock Springer-Verlag, 1992.

\bibitem{M-S}
N.~Manolache and F.-O. Schreyer.
\newblock Moduli of $(1,7)$-polarized abelian surfaces via syzygies.
\newblock {\em Math. Nachr.}, 226:177--203, 2001.

\bibitem{AMa}
A.~Marini.
\newblock {\em On the degenerations of $(1,7)$-polarised abelian surfaces}.
\newblock PhD thesis, Bath University, 2002.

\bibitem{Mie}
M.~Miele.
\newblock {\em Klassifikation der Durchschnitte Heisenberg-Invariante Systeme
  von Quadriken in $\p{}6$}.
\newblock PhD thesis, Erlangen, 1993.

\bibitem{Mu}
S.~Mukai.
\newblock Fano 3-folds.
\newblock In {\em Vectors Bundles and Special Proj. Embeddings (Bergen 1989)},
  volume 179 of {\em Lect. Notes Ser.}, pages 255--263. Lond. Math. Soc., 1992.

\end{thebibliography}
